\date{}
\documentclass{icmart}
\usepackage{psfig,amsmath,amssymb,xspace,url,picinpar,graphicx}
\title[Conformally invariant scaling limits]{Conformally invariant scaling limits\\{\large an overview and a collection of problems}}
\author[Oded Schramm]{Oded Schramm}
\contact{Microsoft Corporation, Redmond, Washington, USA}
\newif\ifhyper\IfFileExists{hyperref.sty}{\hypertrue}{\hyperfalse}
\hyperfalse
\hypertrue
\ifhyper\usepackage{hyperref}\fi
\newcommand{\old}[1]{}

\newtheorem{theorem}{Theorem}[section]

\newtheorem{problem}[theorem]{Problem}
\newtheorem{definition}[theorem]{Definition}

\newcommand{\R}{{\mathbb R}}
\newcommand{\N}{{\mathbb N}}

\newcommand{\C}{{\mathbb C}}
\def\H{\mathbb H}
\def\U{\mathbb U}
\newcommand{\Z}{{\mathbb Z}}

\newcommand{\length}{{\operatorname{length}}}

\newcommand{\eps}{\varepsilon}
\def\eps{\epsilon}
\def\SLE/{SLE}
\def\SLEk#1/{$\mathrm{SLE}_{#1}$}
\def\pa{\partial}
\def\closure#1{\overline{#1}}

\def\TG{\mathrm{TG}}
\def\P{\mathbf{P}}

\def\eref#1{(\ref{#1})}
\def\Bb#1#2{{\def\md{\bigm| }#1\bigl[#2\bigr]}}

\def\Pb{\Bb\P}

\def\bl{\bigl}\def\br{\bigr}

\def\HaraSladeSAW{MR1171762}
\def\DLAdef{MR704464}
\def\DLAsurvey{MR1205417}
\def\MakarovCarleson{MR1815718}
\def\KestenDLA{MR1077203}
\def\LanglandsEtAl{MR1230963}
\def\HPSdynamical{MR1465800}
\def\TsirelsonFourierWalsh{math.PR/9903068}
\def\TsirelsonScalingLimit{math.PR/9903121}
\def\BousquetMelouSchaeffer{math.CO/0211070}
\def\MarckertMokkadem{math.PR/0403398}
\def\ChassaingSchaeffer{MR2031225}
\def\AngleGrowthUIPT{MR2024412}
\def\MadrasSlade{MR1197356}
\def\KennedyVariation{math.PR/0510604}
\def\FourierSurveys{KalaiSafra}
\def\BKSnoise{MR2001m:60016}
\def\KPZ{MR947880}
\def\Ziff{MR1745857}
\def\BollobasRiordanVoronoi{math.PR/0410336}
\def\CFNnearcrit{cond-mat/0510740}
\def\LSWuptoconst{MR1901950}
\def\RSWboolean{MR1398054}
\def\RSWvacant{MR1071809}
\def\WernerLoops{math.PR/0511605}
\def\BenjaminiSchrammVoronoi{MR1646475}
\def\SheffieldGFFsurvey{math.PR/0312099}
\def\KozmaLERWspace{math.PR/0508344}
\def\KozmaLERWplane{math.PR/0212338}
\def\Watts{MR98a:82059}
\def\DubedatWatts{math.PR/0405074}

\def\TaylorSurvey{MR857718}
\def\NienhuisOn{NienhuisOn}
\def\SSexceptional{math.PR/0504586}
\def\FortuinKasteleyn{MR50:12107}
\def\ItzyksonDrouffe{MR1175176}
\def\BeffaraSidoravicius{math.PR/0507220}
\def\BeffaraSLEdim{math.PR/0211322}
\def\KenyonConformalDomino{MR1782431}

\def\MeesterRoy{MR1409145}
\def\GrimmettPercolationEdtwo{MR1707339}
\def\CamiaNewman{math.PR/0504036}
\def\KestenICM{MR1989192}
\def\LevyConfInvBM{MR0029120}

\def\SmirnovPerc{MR1851632}

\def\CardySLESurvey{MR2148644}

\def\LSWsaw{MR2112127}
\def\LSWi{MR2002m:60159a}
\def\LSWii{MR2002m:60159b}
\def\LSWiii{MR1899232}

\def\LSWlesl{math.PR/0112234}
\def\LSWlesl{MR2044671}
\def\LSWrestriction{MR1992830}
\def\LSWonearm{MR2002k:60204}
\def\SWpercexpo{math.PR/0109120}

\def\SchSLE{MR1776084}

\def\CardyFormula{MR92m:82048}

\def\Dudley{MR91g:60001}

\def\RSsle{MR2153402}

\def\LawlerLERWdef{MR81j:60081}
\def\LawlerWernerUniversality{MR2002g:60123}
\def\LawlerWernerBMexponents{MR2000k:60165}
\def\LWi{\LawlerWernerBMexponents}
\def\LWii{\LawlerWernerUniversality}
\def\LawlerSLbook{MR2129588}

\def\PemantleUST{MR92g:60014}

\def\WilsonAlg{MR1427525}
\def\LyonsUSFsurvey{MR99e:60029}
\def\KenyonLERW{MR1819995}
\def\KenyonUST{MR1757962}

\def\LawlerFrontier{MR97g:60110}
\def\SchrammSheffieldHE{MR2184093}
\def\SchrammSheffieldDGFF{DGFF}
\def\SchrammWilsonCoord{math.PR/0505368}
\def\multiplyConnectedSLE{\DZhan,MR2111355,math.PR/0412060,math.PR/0503178}
\def\DZhan{MR2128237}
\def\WernerICM{WernerICM}
\def\LWsoup{MR2045953}
\def\WernerBMclusters{MR2023758}
\def\SladeKestenFest{MR1703123}

\def\BrydgesImbrie{MR2031859}
\def\KenyonDominoGFF{MR1872739}
\def\KestenPc{MR82c:60179}
\def\HarrisPc{MR0115221}
\def\WiermanPc{MR612205}
\def\RussoPercolation{MR0488383}

\def\WernerRestrictionSurvey{MR2178043}
\def\WernerStFlour{MR2079672}

\def\SeymourWelsh{MR58:13410}

\def\NaddafSpencer{MR1461951}
\def\KenyonDominoGFF{MR1872739}

\def\DuplantierKwon{DuplantierKwon}

\def\KagerNienhuisSurvey{MR2065722}

\def\figdir{./}
\def\noopsort#1{}

\begin{document}

\begin{abstract}
Many mathematical models of statistical physics in two dimensions are
either known or conjectured to exhibit conformal invariance.
Over the years, physicists proposed predictions of various
exponents describing the behavior of these models.
Only recently have some of these predictions become accessible to
mathematical proof. One of the new developments is the discovery
of a one-parameter family of random curves
called Stochastic Loewner evolution
or \SLE/. The \SLE/ curves appear as limits of interfaces or paths
occurring in a variety of statistical physics models as the
mesh of the grid on which the model is defined tends to zero.

The main purpose of this article is to list a collection of
open problems.  Some of the open problems indicate aspects of the
physics knowledge that have not yet been understood mathematically.
Other problems are questions about the nature of the \SLE/ curves
themselves. 
Before we present the open problems, the definition of \SLE/
will be motivated and explained, and a brief sketch
of recent results will be presented.
\end{abstract}

\begin{classification}
Primary 60K35; Secondary 82B20, 82B43, 30C35.
\end{classification}

\begin{keywords}
Statistical physics, conformal invariance, stochastic Loewner evolutions, percolation.
\end{keywords}

\maketitle

\section{Introduction}
In the past several years, many predictions from physics regarding the large-scale
behavior of random systems defined on a lattice in two dimensions have become
accessible to mathematical study and proof. Of central importance is the 
asymptotic conformal invariance of these systems. 
It turns out that paths associated with these random configurations
often fall into a one-parameter family of conformally invariant random curves called
stochastic Loewner evolutions, or \SLE/.
We start by motivating \SLE/ through a simple mathematical model of percolation.
After giving the definition of \SLE/, we present a narrative of recent developments.
However, since there are good surveys on the subject
in the literature~\cite{\WernerStFlour,\KagerNienhuisSurvey,\CardySLESurvey,%
\LawlerSLbook,\WernerRestrictionSurvey}, this introductory
part of the paper will be short and cursory.
The rest of the paper will consist of an annotated list of open problems in the subject.

\subsection{Motivation and definition of SLE}

\begin{figure}
\centerline{\psfig{file=\figdir/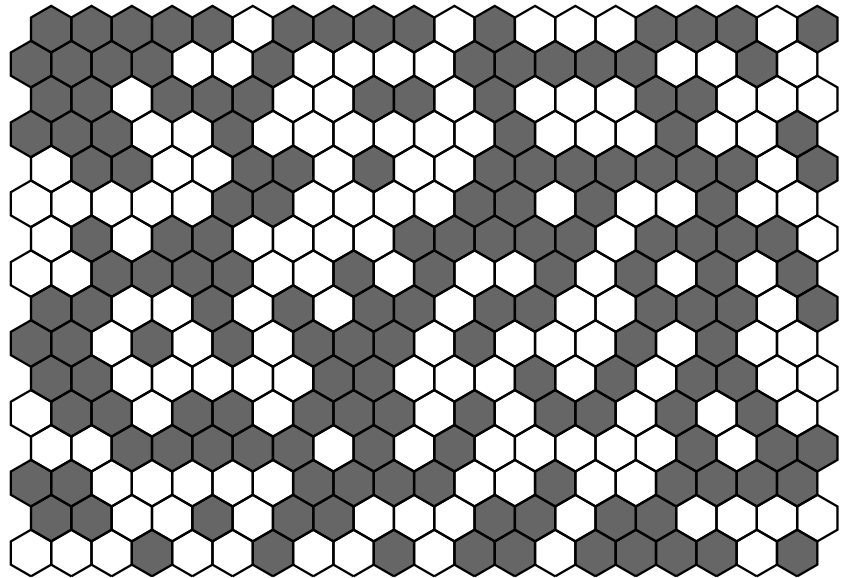,height=3in}}
\caption{\label{f.perc}Site percolation on the triangular grid 
as represented by colored hexagons.}
\end{figure}

To motivate SLE, we now discuss percolation.
More specifically, we define one particular model of percolation
in two dimensions.
Fix a number $p\in[0,1]$. Let $\omega$ be a random subset
of the set of vertices in the triangular grid $\TG$, where 
for vertices $v\in V(\TG)$ the events $v\in\omega$ are independent
and have probability $p$.
In percolation theory one studies the connected
components (a.k.a.\ clusters) of the random subgraph of $\TG$ whose vertex
set is $\omega$ and whose edges are the edges in $\TG$
connecting two elements of $\omega$. Equivalently,
one may study the connected components of the set of white
hexagons in Figure~\ref{f.perc}, where each hexagon
in the hexagonal grid dual to $\TG$ represents
a vertex of $\TG$ and hexagons corresponding to vertices
in $\omega$ are colored white. The reasons for considering this
dual representation are that the figures come out nicer and
that it makes some important definitions more concise.

The above percolation model is site (or vertex) percolation on the
triangular grid. Likewise, there is a bond (or edge) model,
where one considers a random subgraph of a grid whose
vertex set is the set of all vertices of the grid, but
where each edge of the grid is in the percolation subgraph
with probability $p$, independently.
Additionally, there are various percolation models which
are not based on a lattice. Some of these will be discussed
in later sections.

There is an important value $p_c$ of the parameter $p$, which
is the threshold for the existence of an unbounded cluster
and is called the critical value of $p$.
The actual value of $p_c$ varies depending on the particular
percolation model. For site percolation on the triangular
grid, as well as for bond percolation on the square grid,
we have $p_c=1/2$. This is a theorem of Kesten~\cite{\KestenPc},
based on earlier work by Harris~\cite{\HarrisPc}, Russo~\cite{\RussoPercolation}
and Seymour and Welsh~\cite{\SeymourWelsh}. The underlying reason for this
nice value of $p_c$ is a duality which these two models have,
though the precise form of the duality they exhibit is different.
For bond percolation on the triangular grid $p_c=2\,\sin(\pi/18)$~\cite{\WiermanPc}, while
for site percolation on the square grid there is not even a prediction
for the value of $p_c$, though rigorous and experimental estimates exist.
As $p$ increases beyond $p_c$, the large scale behavior of 
percolation undergoes a rapid change.
This is perhaps the mathematically simplest model of a phase transition.
From now on, we will focus our attention on
critical percolation, that is, percolation with $p=p_c$, which is in
many ways the most interesting value of $p$.

\begin{figure}
\centerline{\psfig{file=\figdir/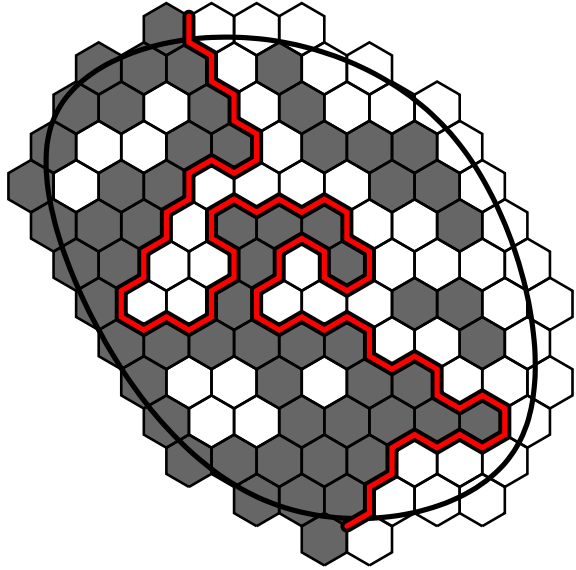,height=3in}}
\caption{\label{f.domperc}The interface associated with percolation.}
\end{figure}

We now define and discuss the percolation interface curve indicated in
Figure~\ref{f.domperc}.
Consider a bounded domain $D$ in the plane $\R^2=\C$ whose boundary is a 
simple closed curve. Let $\pa_+\subset\pa D$ be a proper arc on the boundary of $D$.
Given $\eps>0$ we may consider the collection of hexagons in a hexagonal grid of
mesh $\eps$ which intersect $\closure D$.
Each of these hexagons which meets $\pa_+$ we color white,
each of the hexagons which meet $\pa D$ but not $\pa_+$ we color black,
and each of the hexagons contained in $D$ we color white or black with
probability $1/2$, independently.
In addition to white clusters (connected components of white hexagons)
sometimes, black clusters are also  considered.
Percolation theory is the study of connected components of random sets,
such as these clusters.

For simplicity, we assume that $\pa D$ is sufficiently smooth and
$\eps$ is sufficiently small so that the union of hexagons
intersecting $\pa D$ but not $\pa_+$ is connected.
There is a unique (random) path $\beta$,
which is the common boundary of the white cluster 
meeting $\pa_+$ and the black cluster meeting 
$\pa D$.

The law $\mu_{D,\pa_+,\eps}$ of $\beta$ is a probability measure on the
space of closed subsets of $\closure D$ with the Hausdorff metric.
Smirnov~\cite{\SmirnovPerc} proved that as $\eps\searrow0$ the 
measure $\mu_{D,\pa_+,\eps}$ converges weakly to a measure
$\mu_{D,\pa_+}$, and that $\mu$ is conformally invariant,
in the following sense.
If $f:\closure D\to \closure{D'}\subset\R^2$ is a homeomorphism
that is analytic in $D$, then the push forward of $\mu_{D,\pa_+}$ under
$f$ is $\mu_{f(D),f(\pa_+)}$.
In other words, if $\eps$ is small, then $f(\beta)$ is a good approximation
for the corresponding path defined using a hexagonal grid of mesh $\eps$
in $f(D)=D'$.
This type of conformal invariance was believed to hold for many
\lq\lq critical\rq\rq\ random systems in two dimensions.
However, the only previous result establishing conformal invariance
for a random scaling limit is L\'evy's theorem~\cite{\LevyConfInvBM}
stating that for two-dimensional Brownian motion,
the scaling limit of simple random walk on $\Z^2$, is conformally invariant
up to a time-change.

Smirnov's proof is very beautiful, and the result is
important, but describing the proof will throw
us too far off course (since for this paper percolation is
just an example model, not the primary topic).
The interested reader is encouraged to consult~\cite{\BeffaraSidoravicius,%
\GrimmettPercolationEdtwo,\MeesterRoy,\KestenICM}
for background in percolation and highlights of percolation theory.
An elegant simplification of parts of Smirnov's proof has been
discovered by Vincent Beffara~\cite{BeffaraSmirnovThm}.
A more detailed version of other parts of Smirnov's argument appears
in~\cite{\CamiaNewman}.

\begin{figure}
\centerline{\psfig{file=\figdir/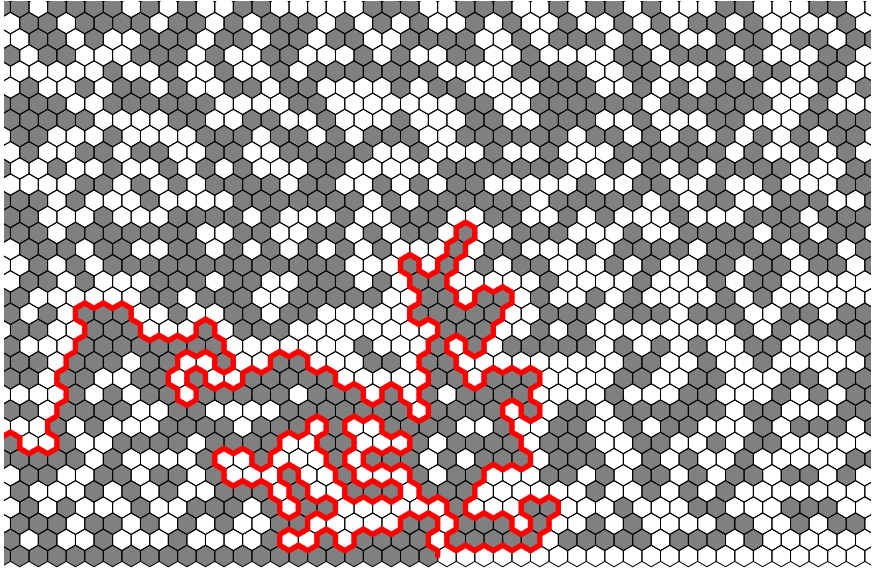,height=3in}}
\caption{\label{f.pcurve}The percolation interface in the upper half plane.}
\end{figure}

Though this was not the original inspiration, 
we will now use Smirnov's result to motivate the definition of SLE.
By conformal invariance, we may venture to understand $\beta$
in the domain of our choice. The simplest situation turns out
to be when $D=\H$ is the upper half plane and $\pa_+$ is the
positive real ray $\R_+$, as in Figure~\ref{f.pcurve}.
(Though this domain $D$ is unbounded, that does not cause any problems.)
We may consider the discrete path
$\beta$ as a simple path $\beta:[0,T)\to\closure \H$ starting near $0$ and
satisfying $\lim_{t\to T}|\beta(t)|=\infty$ (where $T$ is finite or infinite).

We would like to learn about $\beta$ by understanding the 
one-parameter family of conformal maps $g_t$ mapping
$\H\setminus \beta([0,t])$ onto $\H$.
To facilitate this, we must first recall a few basic facts and
discuss Loewner's theorem.
At this point, assume only that $\beta$ is a simple path
in $\closure{\H}$ with $\beta(0)\in\R$
and $\beta(t)\notin\R$ for $t>0$.
The existence of the conformal maps $g_t:\H\setminus\beta([0,t])\to\H$
is guaranteed by Riemann's mapping theorem.
However, $g_t$ is not unique.  In order to choose a specific
$g_t$ for every $t$, we first require that $g_t(\infty)=\infty$.
Schwarz reflection in the real axis implies that $g_t$ is analytic
in a neighborhood of $\infty$, and therefore admits a power series
representation in $1/z$, 
$$
g_t(z)=a_1\,z+a_0+a_{-1}\,z^{-1}+a_{-2}\,z^{-2}+\cdots,
$$
valid for all $z$ sufficiently large. Since $g_t$ maps the real line
near $\infty$ into the real line, it follows that
$a_j\in\R$ for all $j$, and because $g_t:\H\setminus\beta([0,t])\to\H$,
we find that $a_1>0$. We now pick a specific $g_t$ by imposing the
so-called hydrodynamic normalization at $\infty$, namely $a_1=1$ and
$a_0=0$. This can clearly be achieved  by post-composing
with a map of the form $z\mapsto a\,z+b$, $a>0$, $b\in\R$.

The coefficients $a_j$ of the series expansion of $g_t$ are now functions of $t$.
It is not hard to verify that $a_{-1}(t)$ is a continuous, strictly increasing function
of $t$. Clearly, $g_0(z)=z$ and hence $a_{-1}(0)=0$.
We may therefore reparametrize $\beta$ so as to have
$a_{-1}(t)=2\,t$ for all $t>0$. 
This is called the half-plane capacity parametrization of $\beta$.
With this parametrization, a variant of Loewner's theorem~\cite{Loewner} states
that the maps $g_t$ satisfy the differential equation
\begin{equation}
\label{e.chordal}
\frac{dg_t(z)}{dt}= \frac{2}{g_t(z)-W(t)}\,,
\end{equation}
where $W(t):=g_t\bl(\beta(t)\br)$
is called the Loewner driving term. A few comments are in order.
\begin{enumerate}
\item
Although $g_t$ is defined in $\H\setminus \beta([0,t])$,
it does extend continuously to $\beta(t)$, and therefore
$W(t)$ is well defined. 
\item
If $z=\beta(s)$ for some $s$,
then~\eref{e.chordal} makes sense only as long as $t<s$.
That is to be expected. The domain of definition of $g_t$ is shrinking as $t$ increases.
A point $z$ falls out of the domain of $g_t$ at the first time $\tau=\tau_z$ such that
$\liminf_{t\nearrow\tau} g_t(z)-W(t)=0$.
\item The main point here is that information about the path $\beta$ is encoded
 in $W(t)$, which is a path in $\R$.
\item The proof of~\eref{e.chordal} is not too hard. In~\cite[Theorem 2.6]{\LSWi} 
a proof (of a generalization) may be found.
\end{enumerate}
  
\begin{figure}
\centerline{\psfig{file=\figdir/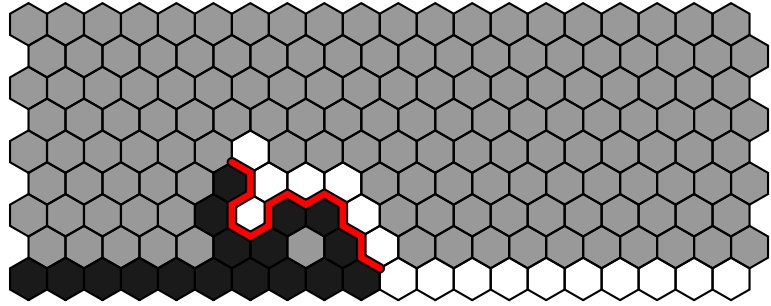,width=4in}}
\caption{\label{f.pcurveini}Initial segment of interface.}
\end{figure}

We now return to the situation where $\beta$ is the percolation interface chosen according
to $\mu_{\H,\R_+,\eps}$, parametrized by half-plane capacity.
It is easy to see that in this case $T=\infty$.
Fix some $s>0$. Suppose that we examine the colors of only those hexagons that
are necessary to determine $\beta([0,s])$. This can be done by sequentially
testing the hexagons adjacent to $\beta$ starting from $\beta(0)$
as follows.  Each time the already determined arc of $\beta$ meets a hexagon whose color
has not yet been examined, we test the color (which permits us to extend the
determined initial arc of $\beta$ by at least one segment), until
$\beta([0,s])$ has been determined. See Figure~\ref{f.pcurveini}.

Now comes the main point. Let $D_s$ be the unbounded component
of the collection of hexagons of undetermined color in $\H$, and let
$\pa_+^s$ be the subset of $\pa D_s$ lying on the boundary
of hexagons of determined color white.
Then the distribution of the continuation $\beta([s,\infty))$ of the interface
given $\beta([0,s])$ is $\mu_{D_s,\pa_+^s,\eps}$.
By Smirnov's theorem, if $G:D_s\to\H$ is the conformal map 
satisfying the hydrodynamic normalization, then the image under
$G$ of $\mu_{D_s,\pa_+^s,\eps}$ is close to $\mu_{\H,G(\pa_+^s),\eps}$.
(Actually, to justify this, one needs a slightly stronger \lq\lq uniform\rq\rq\
version of Smirnov's theorem. But here we want to convey the main ideas, and
do not bother about being entirely precise.)
Now, since $D_s$ approximates $\H\setminus \beta([0,s])$,
it follows that $G$ is very close to $g_s$ and $G(\pa_+^s)$ is close
to $[W(s),\infty)=[g_s(\beta(s)),\infty)$. Therefore, in the limit as $\eps\searrow 0$, we have
for $\beta$ sampled according to $\mu_{\H,\R_+}$ that given
$\beta([0,s])$ the distribution of $ g_s\circ \beta([s,\infty))$
(which is the conformal image of the continuation of the path)
is $\mu_{\H,\R_+}$ translated by $W(s)$.

The Loewner driving term of the path 
$t\mapsto g_s\circ\beta(s+t)$ is $W(s+t)$,
because $g_{s+t}\circ g_s^{-1}$ maps
$\H\setminus g_s\bl(\beta([s,s+t])\br)$ onto $\H$.
The conclusion of the previous paragraph therefore implies that given $(W(t):t\in[0,s])$
the distribution of the continuation of $W$ is identical to the
original distribution of $W$ translated to start at $W(s)$.
This is a very strong property. 
Indeed, for every $n\in\N$ and $t>0$ we
may write $W(t)=\sum_{j=1}^n \bl(W(j\,t/n)-W((j-1)t/n)\br)$,
which by the above is a sum of $n$ independent identically distributed
random variables.
If we assume that the variance of $W(t)$ is finite,
then it is also the sum of the variances of the summands.
By the central limit theorem, $W(t)$ is therefore a Gaussian
random variable.
By symmetry, $W(t)$ has the same distribution as $-W(t)$,
and so $W(t)$ is a centered Gaussian. It now easily follows
that there is some constant $\kappa\ge 0$ such that
$W(t)$ has the same distribution as $B(\kappa\,t)$,
where $B$ is one-dimensional Brownian motion starting at $B(0)=0$.
Using results from the theory of stochastic processes
(e.g., the characterization of continuous martingales as time-changed
Brownian motion), the same conclusion can be reached while replacing
the assumption that $W(t)$ has finite variance with
the continuity of $W(t)$ in $t$.

\bigskip

We have just seen that Smirnov's theorem implies that
the Loewner driving term of a sample from $\mu_{\H,\R_+}$
is $B(\kappa\,t)$ for some $\kappa\ge 0$. This should serve
as adequate motivation for the following definition from~\cite{\SchSLE}.

\begin{definition}
Fix some $\kappa\ge 0$, and let $g_t$ be the solution of
Loewner's equation~\eref{e.chordal}
satisfying $g_0(z)=z$ with $W(t)=B(\kappa\,t)$, where $B$
is standard one-dimensional Brownian motion starting at $B(0)=0$.
Then $(g_t:t\ge 0)$ is called {\bf chordal Stochastic Loewner evolution}
with parameter $\kappa$ or \SLEk\kappa/.
\end{definition}

Of course, \SLEk\kappa/ is a random one-parameter family of maps;
the randomness is entirely due to the Brownian motion.

It has been proven~\cite{\RSsle,\LSWlesl} that with probability $1$
there is a (unique) random continuous path $\gamma(t)$ such that
for each $t\ge 0$ the domain of definition of $g_t$ is
the unbounded component of $\H\setminus\gamma([0,t])$.
The path is given by $\gamma(t)=g_t^{-1}(W(t))$, but
proving that $g_t^{-1}(W(t))$ is well defined is not easy.
It is also known~\cite{\RSsle} that a.s.\ $\gamma(t)$ is a simple
path if and only if $\kappa\le 4$ and is space-filling if and only
if $\kappa\ge 8$.
Sometimes the path $\gamma$ itself is called \SLEk\kappa/.
This is not too inconsistent, because $g_t$ can be reconstructed
from $\gamma([0,t])$ and vice versa.

If $D$ is a simply connected domain in the plane and $a,b\in\pa D$ are
two distinct points (or rather prime ends), then chordal \SLE/ from
$a$ to $b$ in $D$ is defined as the image of $\gamma$ under a conformal
map from $\H$ to $D$ taking $0$ to $a$ and $\infty$ to $b$.
Though the map is not unique, the choice of the map does not
effect the law of the \SLE/ in $D$. This follows from the easily
verified fact that up to a
rescaling of time, the law of the \SLE/ path is invariant
under scaling by a positive real constant, as is the case for Brownian
motion.

The reason for calling the \SLE/ \lq\lq chordal\rq\rq\ is
that it connects two boundary points of a domain $D$.
There is another version of \SLE/, which connects
a boundary point to an interior point, called {\bf radial}
\SLE/. Actually, there are a few other variations, but they
all have similar definitions and analogous properties.

\subsection{A historical narrative}

In this subsection we list some works and discoveries related to \SLE/
and random scaling limits in two dimensions.
The following account is not comprehensive. Some of the topics not covered
here are discussed in Wendelin Werner's~\cite{\WernerICM}
contribution to this ICM proceedings.
We start by very briefly discussing the historical background.

\bigskip

In the survey paper~\cite{\LanglandsEtAl}, Langlands, Pouliot and Saint-Aubin
present a collection of intriguing predictions from statistical physics.
They have discussed these predictions and some simulation
data with Aizenman, which prompted him
to conjecture that the critical percolation
crossing probabilities are asymptotically conformally invariant (see \cite{\LanglandsEtAl}).
This means that the probability $Q_\eps(D,\pa_1,\pa_2)$ that there
exists a critical percolation cluster in a domain $D\subset\C$ connecting two boundary arcs
$\pa_1$ and $\pa_2$ on a lattice with mesh $\eps$ has a limit
$Q(D,\pa_1,\pa_2)$ as $\eps\searrow 0$ and that the limit
is conformally invariant, namely,
$Q(D,\pa_1,\pa_2)=Q\bl(f(D),f(\pa_1),f(\pa_2)\br)$
if $f$ is a homeomorphism from $\closure D$ to $f(\closure {D})$
that is conformal in $D$.
This led John Cardy~\cite{\CardyFormula} to propose his formula (involving hypergeometric
functions) for the asymptotic crossing probability in a rectangle between two opposite edges.
The survey~\cite{\LanglandsEtAl} highlighted these predictions
and the role of the conjectured conformal invariance
in critical percolation, as well as several other statistical
physics models in two dimensions.

Prior to \SLE/ there were attempts to use compositions
of conformal slit mappings and even Loewner's equation
in the study of diffusion limited aggregation 
(DLA).  DLA is a random growth process, which
produces a random fractal and is notoriously
hard to analyse mathematically.
(See~\cite{\DLAdef,\DLAsurvey} for a definition and discussion of DLA.)
Makarov and Carleson~\cite{\MakarovCarleson} used Loewner's equation
to study a much simplified deterministic variant of DLA, which is
not fractal, and Hastings and Levitov~\cite{HastingsLevitov}
have used conformal mapping techniques for a non-rigorous
study of more realistic versions of DLA.
Given that the fractals produced by DLA are not conformally
invariant, it is not too surprising that
it is hard to faithfully model DLA using conformal maps.
Harry Kesten~\cite{\KestenDLA} proved that 
the diameter of the planar DLA cluster after $n$ steps grows asymptotically
no faster than $n^{2/3}$,
and this appears to be essentially the only theorem concerning
two-dimensional DLA, 
though several very simplified variants of DLA have been successfully analysed.

\bigskip

\begin{figure}
\centerline{\hskip3.5in\psfig{file=\figdir/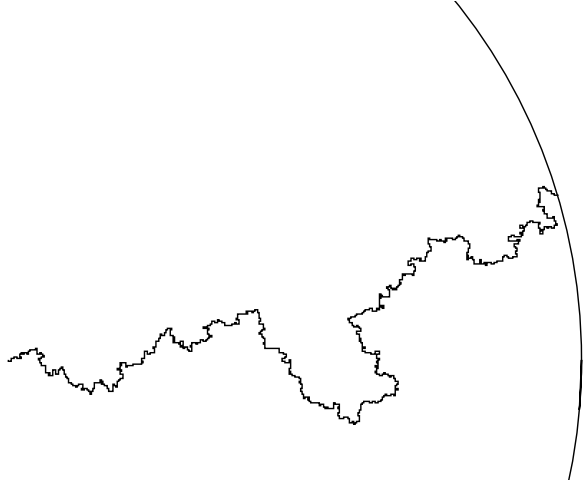,height=2.5in}}
\caption{\label{f.lerw}The LERW in a disk.}
\end{figure}

The original motivation for \SLE/ actually came from
investigating the Loop-erased random walk (a.k.a.\ LERW),
which is a random curve introduced by Greg Lawler~\cite{\LawlerLERWdef}.
Consider some bounded simply-connected domain $D$ in the plane.
Let $G=G(D,\eps)$ be the subgraph of 
a square grid of mesh $\eps$ that falls inside $D$ and let $V_\pa$ be the
set of vertices of $G$ that have fewer than $4$ neighbors in $D$.
Suppose that $0\in D$, and let $o$ be some vertex of $G$ closest to $0$.
Start a simple random walk on $G$ from $o$ (at each step the walk jumps to
any neighbor of the current position with equal probability). We keep
track of the trajectory of the walk at each step, except that every time
a loop is created, it is erased from the trajectory. The walk terminates
when it first reaches $V_\pa$, and the loop-erased random walk from $o$
to $V_\pa$ is the final trajectory.
See Figure~\ref{f.lerw}, where $D$ is a disk.

The LERW is intimately related to the uniform spanning tree.
In particular, if we collapse $V_\pa$ to a single vertex $v_\pa$
and take a random spanning tree of the resulting projection of $G$,
where each possible spanning tree is chosen with equal probability,
then the unique path in the tree joining $o$ to $v_\pa$ (as a set
of edges) has precisely the same law as the LERW
from $o$ to $V_\pa$~\cite{\PemantleUST}.
 This is not a particular property of the square
grid, the corresponding analog holds in an arbitrary finite graph.
In the other direction, there is a marvelous algorithm discovered
by David Wilson~\cite{\WilsonAlg} which builds the uniform spanning
tree by successively adding loop-erased random walks.
The survey~\cite{\LyonsUSFsurvey} is a good window into the
beautiful theory of uniform spanning trees and forests.

Using sophisticated determinant calculations
and Temperley's bijection between the collection of spanning trees and
a certain collection of dimer tilings (which is special to the
planar setting),
Richard Kenyon~\cite{\KenyonLERW,\KenyonUST,\KenyonConformalDomino}
was able to calculate several properties of the LERW.
For example, it was shown that the variance of the winding number
of the above LERW in $G(D,\eps)$ is $(2+o(1))\,\log(1/\eps)$
as $\eps\searrow0$, and that the growth exponent for the number
of edges in a LERW is $5/4$.

\medskip
In~\cite{\SchSLE} it was shown that {\it if\/} the limit of the law of
the LERW as $\eps\searrow0$ exists and is conformally invariant,
then it is a radial \SLEk2/ path, when parametrized
by capacity.
(See subsection~\ref{s.notions} for a description of two alternative
topologies on spaces of probability measures on curves,
for which this convergence may be stated.)
In broad strokes, the reason why it should be
an \SLE/ path is basically the same as the
argument presented above for the percolation interface.
Two important properties of the percolation interface
scaling limit were crucial in the above argument:
conformal invariance and the following Markovian property.
If we condition on an initial segment of the path,
the remainder is an instance of the path in the 
domain slitted by the initial segment starting from
the endpoint  of the initial segment.
Conformal invariance was believed to hold for the LERW
scaling limit, while the Markovian property
does hold for the reversal of the LERW.

The identification of the
correct value of the parameter $\kappa$ as $2$ is based on 
Kenyon's calculated LERW winding variance growth rate and
a calculation of the variance of the winding number
of the radial \SLEk\kappa/ path truncated at distance
$\eps$ from the interior target point.
The latter grows like $(\kappa+o(1))\,\log(1/\eps)$.

It was also conjectured in~\cite{\SchSLE} that the percolation
interface discussed above converges to \SLEk6/.
The identification of the parameter $\kappa$ as $6$ in this case
was based on Cardy's formula~\cite{\CardyFormula}
and the verification that the corresponding
formula holds for \SLEk\kappa/ if and only if $\kappa=6$.

The percolation interface satisfies the following locality property.
The evolution of the path (given its past) does not depend on the
shape of the domain away from the current location of the endpoint
of the path. Though this is essentially obvious, it should be noted
that other interesting paths (such as the LERW scaling limit) do
not satisfy locality. 

Another process that clearly satisfies locality is Brownian motion.
Greg Lawler and Wendelin Werner~\cite{\LWi,\LWii} studied the intersection
exponents of planar Brownian motion and the relations between them.
An example of an intersection exponent is the unique number
$\xi(1,1)$ such that the probability that the paths of two independent
Brownian motions started at distance $1$ apart within the unit
disk $\U$ and stopped when they first hit the circle $R\,\pa\U$
do not intersect one another is $R^{-\xi(1,1)+o(1)}$ as $R\to\infty$.
At the time, there were conjectures~\cite{\DuplantierKwon} for the
values of many of these exponents, which were rational numbers,
but only two of these could be proved rigorously (not accidentally,
those had values $1$ and $2$).
These exponents encode many fundamental properties of Brownian motion.
For example, Lawler~\cite{\LawlerFrontier} showed that the dimension of the outer
boundary of planar Brownian motion stopped at time $t=1$, say,
is $2\,(1-\alpha)$ for a certain intersection exponent $\alpha$.
Lawler and Werner~\cite{\LWi,\LWii} have proved certain relations
between the intersection exponents, and have shown that any process
which like Brownian motion satisfies conformal invariance and
a certain version of the locality property necessarily has
intersection exponents that are very simply related to the 
Brownian exponents.

Since \SLEk6/ was believed to be the scaling limit
of the percolation interface, it should satisfy locality.
It is also conformally invariant by definition.
Thus, the Brownian exponents should apply to \SLEk6/.
Indeed, in a series of papers~\cite{\LSWi,\LSWii,\LSWiii} Lawler,
Werner and the present author proved the conjectured values of
the Brownian exponents by calculating the corresponding exponents
for \SLEk6/ (and using the previous work by Lawler and Werner).
Very roughly, one can say that the reason why the exponents
of \SLE/ are easier to calculate than the Brownian exponents
is that the \SLE/ path, though it may hit itself,
does not cross itself. Thus, the outer boundary of the
\SLE/ path is drawn essentially in chronological order.

Later~\cite{\LSWrestriction} it became clear that the relation between 
\SLEk6/ and Brownian motion is even closer than previously
apparent: the outer boundary of Brownian motion started from
$0$ and stopped on hitting the unit circle $\pa\U$ has the
same distribution as the outer boundary of a variant of
\SLEk6/.

Lennart Carleson observed that, assuming conformal invariance, Cardy's formula 
is equivalent to the statement that $Q(D,\pa_1,\pa_2)=\length(\pa_2)$ when $D$ is an
equilateral triangle of sidelength $1$, $\pa_1$ is its base, and
$\pa_2\subset\pa D$ is a line segment having the vertex
opposite to $\pa_1$ as one of its endpoints.
Smirnov~\cite{\SmirnovPerc} proved Carleson's form of
Cardy's formula for
critical site percolation on the triangular lattice
(that is, the same percolation model we have described
above) and showed that crossing probabilities between 
two arcs on the boundary of a simply connected domain are
asymptotically conformally invariant. 
As a corollary, Smirnov concluded
that the scaling limit of the percolation
interface exists and is equal to the \SLEk6/ path.
This connection enabled proving many conjectures about
this percolation model. For example, the prediction~\cite{DN,NienhuisCoulomb}
that the probability that
the cluster of the origin has diameter larger than $R$ decays like
\begin{equation}\label{e.onearm}
\Pb{\text{the origin is in a cluster of diameter}\ge R}
= R^{-5/48+o(1)}
\end{equation}
as $R\to\infty$ was proved~\cite{\LSWonearm}.
This value $5/48$ is an example of what is commonly referred to
as a critical exponent. Building on earlier work by Kesten
and others, as well as on Smirnov's theorem and \SLE/,
Smirnov and Werner~\cite{\SWpercexpo} were able to determine many useful
percolation exponents.
Julien Dub\'edat~\cite{\DubedatWatts} has used \SLE/ to prove Watts'~\cite{\Watts} formula
for the asymptotic probability that in a given rectangle there are both a
white horizontal and vertical crossing (for the above percolation model
at $p=p_c=1/2$). 

\medskip

\begin{figure}
\centerline{\hfill\psfig{file=\figdir/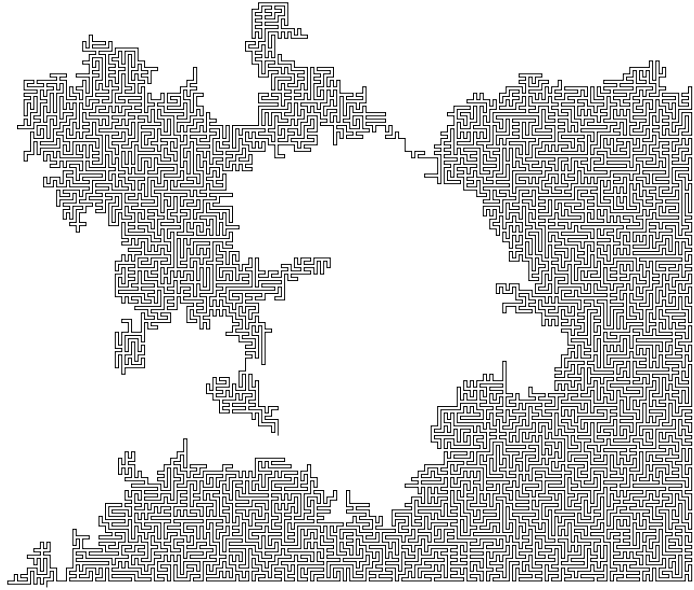,width=\hsize}}
\caption{\label{f.peano}An initial segment of the uniform spanning tree Peano path.}
\end{figure}

The next process for which conformal invariance and
convergence to \SLE/ was established is the LERW~\cite{\LSWlesl}.
Contrary to Smirnov's proof for percolation, where convergence
to \SLE/ was a consequence of conformal invariance, in the case of the
LERW the proof establishes conformal invariance as a consequence
of the convergence to \SLEk2/. More specifically, the argument in~\cite{\LSWlesl}
proceeds by considering the Loewner driving term of the discrete
LERW (before passing to the limit) and proving that the driving term
converges to an appropriately time-scaled Brownian motion.
The same paper also shows that the uniform spanning tree scaling limit is
conformally invariant, and the Peano curve associated with it 
(essentially, the boundary of a thickened uniform spanning tree)
converges to \SLEk8/.
Another difference between the results of~\cite{\SmirnovPerc} and~\cite{\LSWlesl}
is that while the former is restricted to site percolation on the triangular lattice,
the results in~\cite{\LSWlesl} are essentially lattice independent.
Figure~\ref{f.lerw} above shows a fine LERW, which gives an idea of what an
\SLEk2/ looks like. Likewise, Figure~\ref{f.peano} shows a sample
of an initial segment of the uniform spanning tree Peano curve in a rectangular domain.
Note that the curve is space filling, as is \SLEk8/.

Meanwhile, Gady Kozma~\cite{\KozmaLERWplane} came up with a different proof that the
LERW scaling limit exists.
Although Kozma's proof does not identify the limit, it has the advantage
of generalizing to three dimensions~\cite{\KozmaLERWspace}.

\medskip

\begin{figure}
\centerline{\quad\psfig{file=\figdir/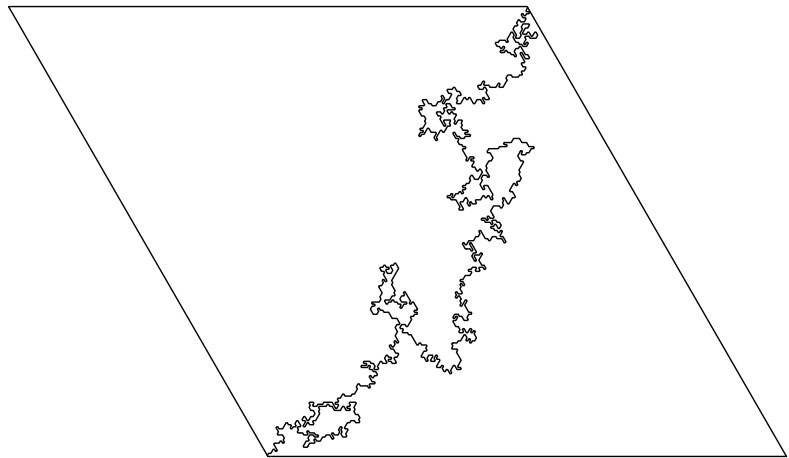,height=2.5in}}
\caption{\label{f.hepath}The harmonic explorer path.}
\end{figure}

There are two discrete models for which convergence to \SLEk4/ has
been established by Scott Sheffield and the present author.
These models are the harmonic explorer~\cite{\SchrammSheffieldHE}
and the interface of the discrete Gaussian free field~\cite{\SchrammSheffieldDGFF}.
The discrete and continuous Gaussian free fields (a.k.a.\ the harmonic crystal)
play an important role in the heuristic
physics analysis of various statistical physics models.
The discrete Gaussian free field is a probability measure
on real valued functions defined on a graph, often a piece of a lattice.
If, for example, the graph is a triangulation of a domain in the plane,
an interface is a curve in the dual graph separating vertices where the function is positive
from vertices where the function is negative.
See~\cite{\SheffieldGFFsurvey} or~\cite{\SchrammSheffieldDGFF} for further details,
and see Figure~\ref{f.hepath} for a simulation of the harmonic explorer, and therefore
an approximation of \SLEk4/.

Sheffield also announced work in progress connecting the Gaussian free field
with \SLEk\kappa/ for other values of $\kappa$. The basic idea is
that while \SLEk4/ may be thought of as a curve solving the equation
$h=0$, where $h$ is the Gaussian free field, for other $\kappa$, the
\SLEk\kappa/ curve may be considered as a solution of 
\begin{equation}\label{e.shef}
c\,\mathrm{winding}(\gamma[0,t])=h(\gamma(t))\,,
\end{equation}
where $c$ is a constant depending on $\kappa$. When $\kappa=4$, the
corresponding constant $c$ is zero, which reduces to the setting
of~\cite{\SchrammSheffieldDGFF}.
Alternatively,~\eref{e.shef} can be heuristically written
as $c\,\gamma'(s)=\exp(i\,h(\gamma(t))$,
where $s$ is the length parameter of $\gamma$. However,
we stress that it is hard to make sense of these equations,
for the Gaussian free field is not a smooth function (in fact,
it is not even a function but rather a distribution).
Likewise, the \SLE/ path is not rectifiable and its winding
at most points is infinite.

\medskip

As mentioned above, there are several different variants of \SLE/ in simply
connected domains:
chordal, radial, as well as a few others, which we have not mentioned.
These variants are rather closely related to one another~\cite{\SchrammWilsonCoord}.
There are also variants defined in the multiply-connected setting~\cite{\multiplyConnectedSLE}.
One motivation for this study comes from statistical physics models, which are easy
to define on multiply connected domains. Since one can easily vary the boundary
conditions on different boundary components of the domain, it is clear that there is
often more than one reasonable choice for the definition of the \SLE/ path.

\medskip

Finally, we mention an intriguing connection between Brownian motion and
\SLEk\kappa/ for $\kappa\in(8/3,4]$. There is the notion of the Brownian
loop soup~\cite{\LWsoup}, which is a Poisson measure on the
space of Brownian motion loops.  
According to~\cite{\WernerBMclusters}, the boundaries of clusters of
a sample from the loop soup measure with intensity $c$ are \SLEk\kappa/-like
paths, where $\kappa=\kappa(c)\in(8/3,4]$.
The proof is to appear in a future joint work of Sheffield and Werner.

\bigskip
The above account describes some of the highlights of the developments
in the field in the past several years.
The rest of the paper will be devoted to a description of
some problems where we hope to see some future progress.
Some of these problems are obvious to anyone working in the field
(though the solution is not obvious), while
others are borrowed from several different sources. A few of the problems
appear here for the first time.
The paper~\cite{\RSsle} contains some additional problems.

\section{Random processes converging to \SLE/}\label{s.sleconv}

As we have seen, paths associated with several random processes have been proved to 
converge to various \SLE/ paths.
However, the list of processes where the convergence is expected but not proved
yet is longer. This section will present questions of this sort, most
of which have previously appeared in the literature.

\medskip
The strategy of the proofs of convergence to \SLE/ in the
papers~\cite{\LSWlesl,\SchrammSheffieldHE,\SchrammSheffieldDGFF} is very similar.
In these papers, a collection of martingales with respect to the filtration given by the
evolution of the curve is used to gain information about the Loewner driving
term of the discrete curve. Although such a proof is also possible for
the percolation interface (using Cardy's formula),
this technique was not available at the time,
and Smirnov used instead an argument which uses the independence properties
of percolation in an essential way and is therefore not likely to be applicable
to many other models. Thus, it seems that presently the most promising technique
is the martingale technique from~\cite{\LSWlesl}.

\subsection{Notions of convergence}\label{s.notions}

To be precise, we must describe the meaning of these scaling limits.
In fact, there are at least two distinct reasonable notions of convergence, which
we now describe.
Suppose that $\gamma_n$ are random paths in the closed upper half plane
$\closure\H$ starting from $0$. Consider the one-point compactification
$\hat\C=\C\cup\{\infty\}$ of $\C=\R^2$, which may be thought of as the sphere $S^2$.
The law of $\gamma_n$ may be thought of as a Borel probability measure
on the Hausdorff space of closed nonempty subsets of $\hat\C$.
Since that Hausdorff space is compact, the space of Borel probability measures on
it is compact with respect to weak convergence of measures~\cite{\Dudley}.
We may say that $\gamma_n$ converges in the Hausdorff sense to a random set
$\gamma\subset\closure\H\cup\{\infty\}$ if the law of $\gamma_n$ converges
weakly to the law of $\gamma$. A similar definition applies to curves in the
closure of a bounded domain $D\subset\C$.
The above stated instances of convergence to \SLE/ hold with respect to this
notion. However, in the case of convergence to \SLEk8/, this does not mean very much,
for \SLEk8/ fills up the domain.

The second notion of convergence is stronger. Suppose that each $\gamma_n$ is a.s.\
continuous with respect to the half-plane capacity parametrization from
$\infty$. (This holds, in particular, if $\gamma_n$ is a.s.\ a [continuous]
simple path.) 
If $d$ is a metric on $\hat \C=\C\cup\{\infty\}$
compatible with its topology, then we may consider
the metric 
$$
d^*(\beta_1,\beta_2):= \sup_{t\in[0,\infty)}d\bl(\beta_1(t),\beta_2(t)\br)
$$
on the space of continuous paths defined on $[0,\infty)$.
We may say that $\gamma_n$ converges to a random path $\gamma$
weakly-uniformly if the law of $\gamma_n$ converges weakly to the law of
$\gamma$ in the space of Borel measures with respect to the metric $d^*$.
This implies Hausdorff convergence. Since $d^*$ is finer than the
Hausdorff metric, there are more functions from the space of paths
to $\R$ that are continuous with respect to $d^*$ than with respect to
the Hausdorff metric. Consequently, weakly-uniform convergence
is stronger than Hausdorff convergence.
In all the results stated above saying that some random path converges to \SLE/,
the convergence is weakly-uniform when the paths are parametrized by capacity
or half-plane capacity (depending on whether the convergence is to radial or
chordal \SLE/, respectively).

In the following, when we ask for convergence to \SLE/, we will mean weakly-uniform
convergence. However, weaker nontrivial forms of convergence would also be very interesting.

\subsection{Self avoiding walk}

Let $G$ be either the square, the hexagonal or the triangular
grid in the plane, positioned so that $0$ is some vertex in $G$.
For $n\in\N$ consider the uniform measures on all
self avoiding $n$-step walks in $G$
that start at $0$ and stay in the upper half plane.
It has been shown in~\cite{\LSWsaw} that when 
$G=\Z^2$ the limiting measure as $n\to\infty$ exists.
(The same proof probably applies for the other alternatives for $G$,
provided that $G$ is positioned so that horizontal lines
through vertices in $G$ do not intersect the relative interior
of edges of $G$ which they do not contain.)

\begin{problem}[\cite{\LSWsaw}]\label{p.saw}
Let $\gamma$ be a sample from the $n\to\infty$ limit of the uniform
measure on $n$-step self avoiding paths in the upper half plane described
above. Prove that the limit as $s\searrow0$ of the law of $s\,\gamma$
exists and that it is \SLEk{8/3}/.
\end{problem}

The convergence may be considered with respect to either of the two
topologies discussed in subsection~\ref{s.notions}.

In~\cite{\LSWsaw} some consequences of this convergence are indicated,
as well as support for the conjecture.

There are some indications that the
setting of the hexagonal lattice is easier: the rate of growth
of the number of self avoiding paths on the hexagonal grid is
predicted~\cite{\NienhuisOn} to be $\bl(2+\sqrt 2+o(1)\br)^{n/2}$;
no such prediction exists
for the square grid or triangular lattice.

In dimensions $d>4$ Takashi Hara and Gordon Slade~\cite{\HaraSladeSAW} proved
that the scaling limit of self-avoiding random
walk is Brownian motion. This is also believed to be the case
for $d=4$. See~\cite{\MadrasSlade} for references and further background
on the self-avoiding walk.

\subsection{Height models}\label{ss.heights}

There is a vast collection of height model
interfaces that should converge to \SLEk4/. 
The one theorem in this regard is the convergence of the interface of the
Gaussian free field~\cite{\SchrammSheffieldDGFF}.
This was motivated by Kenyon's theorem stating that
the domino tiling height function converges to the
Gaussian free field~\cite{\KenyonDominoGFF} and by Kenyon's conjecture
that the double domino interface converges to \SLEk4/ (see~\cite{\RSsle} for a statement
of this problem).

The domino height function is a function on $\Z^2$ associated
with a domino tiling (see~\cite{\KenyonDominoGFF}).
Its distribution is roughly (ignoring boundary issues)
the uniform measure on functions $h:\Z^2\to \Z$
such that $h(0,0)=0$,
$$
h(x,y)\mod 4 =\begin{cases}
0 &\text{$x,y$ even}\\
1 &\text{$x$ odd $y$ even}\\
2 &\text{$x,y$ odd}\\
3 &\text{$x$ even $y$ odd}
\end{cases}
$$
and $|h(z)-h(z')|\in\{1,3\}$ if $|z-z'|=1$, $z,z'\in\Z^2$.

Let $D$ be a bounded domain in the plane whose boundary
is a simple path in the triangular lattice, say.
Let $\pa_+$ and $\pa_-$ be complementary arcs
in $\pa D$ such that  the two common endpoints of these arcs
are midpoints of edges. Consider the uniform measure on functions
$h$ taking odd integer values on vertices in $\closure D$
such that $h=1$ on $\pa_+$, $h=-1$ on $\pa_-$, and
$|h(v)-h(u)|\in\{0,2\}$ for neighbors $v,u$.
We may extend such a function $h$ to $\closure D$ by affine
interpolation within each triangle, and this interpolation
is consistent along the edges. There is then a unique
connected path $\gamma$ that is the connected component
of $h^{-1}(0)$ that contains the two endpoints of each
of the two arcs
$\pa_\pm$.

\begin{problem}\label{p.h1}
Is it true that the path $\gamma$ tends to \SLEk4/?
Does the law of $h$ converge to the Gaussian free field?
\end{problem}

The convergence we expect for $h$ is in the same sense
as in~\cite{\KenyonDominoGFF}.

Note that if we restrict in the above the image of
$h$ to be $\{1,-1\}$, we obtain critical site
percolation on the triangular grid, and the limit
of the corresponding interface is in this case \SLEk6/.

Now suppose that $(D,\pa_+,\pa_-)$ is as above. Let $\lambda\in(0,1/2]$ be
some constant and 
consider the uniform measure on functions
$h$ taking real values on vertices in $\closure D$
such that $h=\lambda$ on $\pa_+$, $h=-\lambda$ on $\pa_-$
and $|h(v)-h(u)|\le 1$ for every edge $[v,u]$.

\begin{problem}
Is it true that for some value of $\lambda$ the corresponding
interface converges to \SLEk4/?
Does the law of $h$ converge in some sense to the
Gaussian free field?
\end{problem}

In the case of the corresponding questions for the
discrete Gaussian free field, there is
just one constant $\lambda$ such that the interface
converges to \SLEk4/. For other choices
of $\lambda$ the interface converges to a
well-known variant of \SLEk4/~\cite{\SchrammSheffieldDGFF}.

There are some restricted classes of height models for which convergence
to the Gaussian free field is known~\cite{\NaddafSpencer}.
It may still be very hard to prove that the corresponding interface
converges to \SLEk4/. One interesting problem of this sort 
is the following.

\begin{problem}[\cite{\SchrammSheffieldDGFF}]
If we project the Gaussian free field 
onto the subspace spanned by the eigenfunctions of the Dirichlet Laplacian
with eigenvalues in $[-r,r]$ and add the harmonic function with boundary
values $\pm\lambda$ on $\pa_\pm$, does the corresponding interface converge to
\SLEk4/ as $r\to\infty$ when $\lambda$ is chosen appropriately?
\end{problem}

The problem is natural, because the Gaussian free field is related to the
Dirichlet Laplacian. In particular, the projections of the field onto the spaces spanned
by eigenfunctions with eigenvalues in two disjoint intervals are independent.

\subsection{The Ising, FK, and $O(n)$ loop models}

The Ising model is a fundamental physics model for magnetism.
Consider again a domain $D$ adapted to the triangular lattice
and a partition $\pa D=\pa_+\cup\pa_-$ as in subsection~\ref{ss.heights}.
Now consider a function $h$ that take the values $\pm1$
on vertices in $\closure D$ such that $h$ is $1$ on $\pa_+$ and $-1$ on $\pa_-$.
On the collection of all such functions we put a probability measure
such that the probability for a given $h$ is
proportional to $e^{-2\beta k}$, where $\beta$ is a parameter and $k$ is the number of edges
$[v,u]$ such that $h(v)\ne h(u)$.
This is known as the Ising model and the value associate to a vertex is
often called a spin.
It is known that the critical value $\beta_c$ (which we do not define here in the context of the Ising model)
for $\beta$ satisfies $e^{2\beta}=\sqrt 3$ (see~\cite{McCoyWu,\ItzyksonDrouffe}).  
Again, the interface at the critical $\beta=\beta_c$ is believed to converge to an
\SLE/ path, this time \SLEk3/. 
For $\beta\in[0,\beta_c)$ fixed, the interface should converge
to \SLEk6/. Note that when $\beta=0$, the model is again identical
to critical site percolation on the triangular grid,
and the interface does converge to \SLEk6/.

\begin{problem}
Prove that when $\beta=\beta_c$, the interface converges to \SLEk3/
and when $\beta\in(0,\beta_c)$ to \SLEk6/.
\end{problem}

When $\beta>\beta_c$ we do not expect convergence to \SLE/,
and do not expect conformal invariance. 
The interface scaling limit is in this case a straight line segment
if the domain is convex~\cite{MR1704276} (see also~\cite{MR2165255}).

\bigskip

The Fortuin-Kasteleyn~\cite{\FortuinKasteleyn} (FK) model
(a.k.a.\ the random cluster model)
is a probability measure on the collection of all subsets of the
set of edges $E$ of a finite graph $G=(V,E)$.
In the FK model, the measure of each $\omega\subset E$ is proportional to
$\bl(p/(1-p)\br)^{|\omega|}\,q^c$, where $q>0$ and $p\in(0,1)$ are parameters,
$|\omega|$ is the
cardinality of $\omega$, and $c$ is the number of connected components
of the subgraph $(V,\omega)$. The FK model is very closely related
to the well known Potts model~\cite{BaxterKellandWu},
which is a generalization of the Ising model.
Many questions about  the Potts model 
can be translated to questions about the FK model and vice versa.

On the grid $\Z^2$, when
$p=\sqrt q/(1+\sqrt q)$ the FK model satisfies a form of self-duality.

\begin{problem}[\cite{\RSsle}]
Prove that
when $q\in(0,4)$ and $p=\sqrt q/(1+\sqrt q)$,
the interface of the FK model on $\Z^2$ with appropriate boundary conditions
converges to \SLEk \kappa/, where
$\kappa=4\,\pi/\cos^{-1}\bl(-\sqrt q/2\br)$.
(See~\cite{\RSsle} for further details.)
\end{problem}

\bigskip

The $O(n)$ loop model on a finite graph $G=(V,E)$ is a measure on the collection
of subgraphs of $G$ where the degree of every vertex in the subgraph is $2$.
(The subgraph does not need to contain all the vertices.)
The probability of each such subgraph is proportional to $x^e\,n^c$,
where $c$ is the number of connected components, $e$ is the number of edges
in the subgraph, and $x,n>0$ are parameters.
When $n$ is a positive integer, the $O(n)$ loop model is derived from the $O(n)$ spin model,
which is a measure on the set of functions which associate to every
vertex a unit vector in $\R^n$. 

In order to pin down a specific long path in the $O(n)$ loop model,
we pick two points on the boundary
of the domain and require that in the random subgraph the degrees of two boundary vertices
near these two points be $1$, while setting the degrees of all other boundary vertices to $0$, say.
Then the measure is supported on configurations with one simple path and a collection of loops.

Now we specialize to the hexagonal lattice.
Set $x_c(n):=\bl(2+\sqrt{2-n}\br)^{-1/2}$,
which is the conjectured critical parameter~\cite{\NienhuisOn}.

\begin{problem}[\cite{\KagerNienhuisSurvey}]
Prove that when $n\in[0,2]$ and $x=x_c(n)$
[respectively, $x>x_c(n)$]
the scaling limit of the path containing the two special boundary vertices
is chordal \SLEk\kappa/, where
$\kappa\in[8/3,4]$ [respectively, $\kappa\in [4,8]$] and $n=-2\,\cos(4\,\pi/\kappa)$.
\end{problem}

When $x<x_c(n)$, we expect the scaling limit to be a straight line segment
(if the domain is convex).
The fact that at $n=1$ we get the same limits as for the Ising model is
no accident. It is not hard to see that the $O(1)$ loop measure coincides with
the law of boundaries of Ising clusters.
See~\cite{\KagerNienhuisSurvey} for further details.

Similar conjectures should hold in other lattices. However, the values of the
critical parameters are expected to be different.

\subsection{Lattice trees}\label{ss.lattree}

We now present an example of a discrete model where we
suspect that perhaps conformal invariance might hold.
However, we do not presently have a candidate for the
scaling limit.

Fix $n\in\N_+$, and consider the collection 
of all trees contained in the grid $G$ that contain
the origin and have $n$ vertices. Select a tree $T$ from
this measure, uniformly at random. 

\begin{problem}
What is the growth rate of the expected diameter
of such a tree?
If we rescale the tree so that the expected (or median)
diameter is $1$, is there a limit for the law of the
tree as $n\to\infty$? What are its geometric and topological
properties? Can the limit be determined?
\end{problem}

It would be good to be able to produce some pictures.
However, we presently do not know how to sample from
this measure.

\begin{problem}
Produce an efficient algorithm which samples lattice
trees approximately uniformly, or prove that such an algorithm
does not exist.
\end{problem}

See~\cite{\SladeKestenFest} for background on lattice trees
and for results in high dimensions and~\cite{\BrydgesImbrie}
for an analysis of a related continuum model.

\subsection{Percolation interface}

It is natural to try to extend the understanding of percolation
at $p_c$ to percolation at a parameter $p$ tending to $p_c$.
One possible framework is as follows.
Fix a parameter $q\in(0,1)$.
Suppose that in Figure~\ref{f.pcurve} with small mesh
$\epsilon>0$ we choose $p(\epsilon)$
so that the probability to have a left to right
crossing of white hexagons in some fixed $1\times 1$
square in the upper half plane is $q$ at percolation parameter $p=p(\epsilon)$.
The corresponding interface will still be an unbounded path starting
at $0$, but its distribution will be different
from the interface at $p=1/2$ if $q\ne 1/2$.
Thus, it is natural to ask

\begin{problem}[Lincoln Chayes (personal communication)]
What is the scaling limit of the interface as $\eps \searrow0$
and $p=p(\eps)$ if $q\ne 1/2$ is fixed?
\end{problem}
This problem is also very closely related to a problem
formulated by F.~Camia, L. Fontes and C.~Newman~\cite{\CFNnearcrit}.

Site percolation on the triangular lattice is only one of
several different models for percolation in the plane.
Among discrete models, widely studied is bond percolation
on the square grid. As for site percolation on the triangular
lattice, the critical probability is again $p_c=1/2$.

At present, Smirnov's proof does not work for bond percolation
on the square grid. The proof uses the invariance of the
model under rotation by $2\,\pi/3$.
 Thus, the following
problem presents itself.

\begin{problem}
Prove Smirnov's theorem for critical bond percolation on $\Z^2$.
\end{problem}

Some progress on this problem has been reported by Vincent Beffara~\cite{BeffaraTalk}.
\medskip

There are other natural percolation models which have been
studied. Among them we mention Voronoi percolation and
the boolean model. In Voronoi percolation one has two
independent Poisson point processes in the plane, $W$ and $B$,
with intensities $p$ and $1-p$, respectively.
Let $\hat W$ be the closure of the set of points in $\R^2$ closer to 
$W$ than to $B$, and let
$\hat B$ be the closure of the set of points closer to $B$.
The set $\hat W$ is a sample from Voronoi percolation
at parameter $p$.
Some form of conformal invariance
was proved for Voronoi percolation~\cite{\BenjaminiSchrammVoronoi},
but the version proved does not imply convergence
to \SLE/. It is neither stronger nor weaker than
the conformal invariance proved by Smirnov.
Notable recent progress has been made for
Voronoi percolation by Bollob\'as and Riordan~\cite{\BollobasRiordanVoronoi},
who established the very useful Russo-Seymour-Welsh
theorem, as well as $p_c=1/2$ for Voronoi percolation.

\begin{problem}
Prove Smirnov's theorem for Voronoi percolation.
\end{problem}

The boolean percolation model (a.k.a.\ continuum percolation)
can be defined by taking a Poisson set
of points $W\subset\R^2$ of intensity $1$ and letting
$\hat W$ be the set of points in the
plane at distance at most $r$ from $W$.
Here, $r$ is the parameter of the model.
(Alternatively, one may fix $r=1$, say, and
let the intensity of the Poisson process be the parameter,
but this is essentially the same, by scaling.)
The Russo-Seymour-Welsh theorem is known for this model~\cite{\RSWboolean,\RSWvacant},
but the critical value of the parameter has not been identified.
A nice feature which the model shares with Voronoi percolation
is invariance under rotations.

\begin{problem}
Prove Smirnov's theorem for boolean percolation.
\end{problem}

\section{Critical exponents}

The determination of critical exponents has been one motivation
to prove conformal invariance for discrete models.
For example, den Nijs and Nienhuis predicted~\cite{DN,NienhuisCoulomb} that the
probability that the critical percolation cluster of the origin 
has diameter larger than $R$ is $R^{-5/48+o(1)}$ as
$R\to\infty$.
Likewise, the probability that a given 
site in the square $[-R,R]^2$ is pivotal
for a left-right crossing of the square $[-2\,R,2\,R]^2$
was predicted to be $R^{-5/4+o(1)}$.
(Here, pivotal means that the occurence or non-occurence of a crossing
would be modified by flipping the status of the site.)
These and other exponents were proved for site percolation on the
triangular grid using Smirnov's theorem and
\SLE/~\cite{\LSWonearm,\SWpercexpo}.
The determination of the exponents is very useful for the
study of percolation.

Richard Kenyon~\cite{\KenyonLERW} calculated by 
enumeration techniques involving determinants
the asymptotics of the probability that an edge
belongs to a loop-erased random walk.
The probability decays like $R^{-3/4+o(1)}$
when the distance from the edge to the endpoints of the
walk is $R$. However, Kenyon's estimate is much more precise;
he shows that, in a specific domain, $R^{3/4}$ times the probability is bounded
away from zero and infinity, and in fact estimates the
probability as $f\,R^{-3/4}\,\bl(1+o(1)\br)$ as $R\to\infty$,
where $f$ is an explicit function of the positiong of the edge.

Thus, it is natural to ask for such precise estimates for
the important percolation events as well. Namely,

\begin{problem}
Improve the estimates $R^{-5/48+o(1)}$ and
$R^{-5/4+o(1)}$ mentioned above (as well as other similar estimates)
to more precise formulas. It would be especially nice to obtain
estimates that are sharp up to multiplicative constants.
\end{problem}

In addition to the case of the loop-erased random walk mentioned
above, estimates up to constants are known for events involving
Brownian motions~\cite{\LSWuptoconst}.
\medskip

The difficulty in getting more precise estimates is not in the
analysis of \SLE/. Rather, it is due to the passage between the
discrete and continuous setting. Consequently, the above
problem seems to be related to the following.

\begin{problem}
Obtain reasonable estimates for the speed of convergence of
the discrete processes which are known to converge to \SLE/.
\end{problem}

There are still critical exponents which do not seem accessible
via an \SLE/ analysis. For example, we may ask

\begin{problem}
Calculate the number $\alpha$ such that on the
event that there is a left-right crossing
in critical percolation in the square $[0,R]^2$,
the expected length of the shortest crossing is
$R^{\alpha+o(1)}$.
\end{problem}

Ziff~\cite{\Ziff} predicts an exponent which is related
to this $\alpha$, but it seems that there is currently
no prediction for the exact value of $\alpha$.

\section{Quantum gravity}

Consider the uniform measure $\mu_n$ on equivalence classes
of $n$-vertex triangulations of the sphere, where
two triangulations are considered equivalent
if there is a homeomorphism of the sphere taking one to the
other. One may view a sample from this measure with
the graph metric as a random geometry on the sphere.
Such models go under the name \lq\lq quantum gravity\rq\rq\ in
physics circles.
One may also impose statistical physics models on such
random triangulations. For example, the sample space
may include such a triangulation (or rather, equivalence
class of triangulations) together with a map $h$ from the
vertices to $\{-1,1\}$. The measure of such a pair
may be taken proportional to $a^k$, where $k$ is the
number of edges $[v,u]$ in the triangulation for
which $h(v)\ne h(u)$ and $a>0$ is a parameter.
Thus, we are in effect considering a triangulation weighted
by the Ising model partition function.
(The partition function is in this case the sum of all the weights
of such functions $h$ on the given triangulation.)
Likewise, one may weight the triangulation by other kinds
of partition functions.

In some cases it is easier to make a heuristic analysis
of such statistical physics models in the quantum gravity
world than in the plane. The enigmatic KPZ formula
of Knizhnik, Polyakov and Zamolodchikov~\cite{\KPZ}
was used in physics to predict properties
of statistical physics in the plane from the corresponding
properties in quantum gravity. Basically, the KPZ formula
is a formula relating exponents in quantum gravity to the
corresponding exponents in plane geometry.

To date, there has been progress in the mathematical
(as well as physical) understanding of the statistical
physics in the plane as well as in quantum gravity~\cite{\AngleGrowthUIPT,math.PR/0501006,\BousquetMelouSchaeffer}.
However, there is still no mathematical understanding of
the KPZ formula. In fact, the author's understanding of KPZ
is too weak to even state a concrete problem.

However, we may ask about the scaling limit of $\mu_n$.
There has been significant progress lately describing some
aspects of the geometry of samples from $\mu_n$~\cite{\AngleGrowthUIPT,\ChassaingSchaeffer}.
In particular, it has been shown by Chassaing and Schaeffer~\cite{\ChassaingSchaeffer}
that if $D_n$ is the graph-metric
diameter of a sample from $\mu_n$, then $D_n/{n^{-1/4}}$
converges in law to some random variable in $(0,\infty)$.
However, the scaling limit of samples from $\mu_n$ is not known.
On the collection of compact metric spaces, we may consider the
Gromov-Hausdorff distance $d_{GH}(X,Y)$, which is the infimum
of the Hausdorff distance between subsets $X^*$ and $Y^*$
in a metric space $Z^*$ over all possible
triples $(X^*,Y^*,Z^*)$
such that $Z^*$ is a metric space,
$X^*,Y^*\subset Z^*$, $X$ is isometric with $X^*$
and $Y$ is isometric with $Y^*$.
Let $X_n$ be a sample from $\mu_n$, considered as a metric space
with the graph metric scaled by $n^{-1/4}$, and let $\mu^*_n$ denote
the law of $X_n$.

\begin{problem}\label{p.UTL}
Show that the weak limit $\lim_{n\to\infty}\mu^*_n$
with respect to the Gromov-Hausdorff metric exists.
Determine the properties of the limit.
\end{problem}

Note that~\cite{\MarckertMokkadem} proves convergence of
samples from $\mu_n$, but the metric used there on
the samples from $\mu_n$ is very different and consequently a
solution of Problem~\ref{p.UTL} does not seem to follow.

\section{Noise sensitivity, Fourier spectrum, and dynamical percolation}

The indicator function of the event of having a percolation crossing
in a domain between two arcs on the boundary is a boolean function of
boolean variables.
Some fundamental results concerning percolation are based on general 
theorems about boolean functions.
(One can mention here the BK inequality,
the Harris-FKG inequality and the Russo formula. See, e.g.~\cite{\GrimmettPercolationEdtwo}.)
Central to the theory of boolean functions is the
Fourier expansion. Basically, if $f:\{-1,1\}^n\to\R$
is any function of $n$ bits, the Fourier-Walsh expansion of $f$ is
$$
f(x)=\sum_{S\subset[n]} \hat f(S)\,\chi_S(x)\,,
$$
where $\chi_S(x)=\prod_{j\in S} x_j$ for $S\subset [n]=\{1,\dots,n\}$.
When we consider $\{-1,1\}^n$ with the uniform probability measure,
the collection $\{\chi_S:S\subset[n]\}$ forms an orthonormal
basis for $L^2(\{-1,1\}^n)$. (Often other measures are also considered.)
If we suppose that $\|f\|_2=1$ (in particular, this holds if
$f:\{-1,1\}^2\to\{-1,1\}$),
then the Parseval identity gives
$\sum_{S\subset[n]} \hat f(S)^2=1$. Thus, we get a probability measure $\mu_f$
on $2^{[n]}=\{S:S\subset[n]\}$ for which $\mu_f(\{S\})=\hat f(S)^2$.
The map $S\mapsto |S|$, assigning to each $S\subset[n]$ its cardinality
pushes forward the measure $\mu_f$ to a measure $\tilde\mu_f$
on $\{0,1,\dots,n\}$. This measure $\tilde \mu_f$ may be called
the Fourier spectrum of $f$. The Fourier spectrum encodes important information
about $f$, and quite a bit of research on the subject
exists~\cite{\FourierSurveys}. For instance, one can read off from $\tilde\mu_f$ the sensitivity
of $f$ to noise (see~\cite{\BKSnoise}).

When $f$ is a percolation crossing function (i.e.,\ $1$ if there is a crossing,
$-1$ otherwise), the corresponding index set $[n]$ is identified
with the collection of relevant sites or bonds, depending if it is a site or bond
model. Though there is some partial understanding of the Fourier spectrum
of percolation~\cite{\SSexceptional}, the complete picture is unclear.
For example, if the domain is approximately an $\ell \times \ell$
square in the triangular lattice, then for the indicator function $f_{\ell}$
of  a crossing in critical
site percolation it is known~\cite{\SSexceptional} that
for every $\alpha<1/8$
\begin{equation}\label{e.percnoise}
\tilde \mu_{f_\ell}([1,\ell^\alpha])\to 0
\end{equation}
as $\ell\to\infty$.
This is proved using the critical exponents for site percolation
as well as 
an estimate for the Fourier coefficients
of general functions (based on the existence of an algorithm
computing the function which is unlikely to examine any specific
input variable).
On the other hand, using the percolation exponents one can show that~\eref{e.percnoise}
fails if $\alpha>3/4$. This is based on calculating the expected number of sites pivotal
for a crossing (i.e., a change of the value of the corresponding input variable would change
the value of the function) as well as showing that the second moment is bounded
by a constant times the square of the expectation. The expected number of pivotals
is known~\cite{\SWpercexpo} to be $\ell^{3/4+o(1)}$ as $\ell\to\infty$.
It is reasonable to conjecture that~\eref{e.percnoise}
holds for every $\alpha<3/4$.

\begin{problem}\label{p.muell}
Is it true that $\lim_{\ell\to\infty}
\tilde\mu_{f_\ell}([\ell^{\alpha_1},\ell^{\alpha_2}])=1$
if $\alpha_1<3/4<\alpha_2$?
Determine the asymptotic behavior of
$\tilde\mu_{f_\ell}([\ell^{\alpha_1},\ell^{\alpha_2}])$ as $\ell\to\infty$
for arbitrary $0\le \alpha_1<\alpha_2\le 2$.
\end{problem}

Estimates on the Fourier coefficients of percolation crossings (in an annulus)
play a central part in the proof~\cite{\SSexceptional} that dynamical percolation
has exceptional times. Dynamical percolation (introduced in~\cite{\HPSdynamical})
is a model in which at each fixed time one sees an ordinary percolation configuration,
but the random bits determining whether or not a site (or bond) is open undergo
random independent flips at a uniform rate, according to
independent Poisson processes. The main result of~\cite{\SSexceptional}
is that dynamical critical site percolation on the triangular lattice has
exceptional times at which there is an infinite percolation component. These set of times
are necessarily of zero Lebesgue measure. 
A better understanding of the Fourier coefficients may lead to sharper
results about dynamical percolation, such as the determination
of the dimension of exceptional times. Some upper and lower bounds for
the dimension are known~\cite{\SSexceptional}.

One would hope to understand the measure $\mu_f$ geometrically.
Gil Kalai (personal communication) has suggested the
problem of determining the scaling limit of $\mu_f$.
More specifically, the Fourier index set $S$ for critical
percolation crossing of a square
is naturally identified with a subset
of the plane. If we rescale the square to have edge length $1$
while refining the mesh, then $\mu_f$ may be thought of as
a probability measure on the Hausdorff space $\mathcal H$ of closed subsets
of the square. It is reasonable to expect that $\mu_f$
converges weakly to some probability measure $\mu$
on $\mathcal H$.
We really do not know what samples from $\mu$ look like.
Could it be that $\mu$ is supported on singletons?
Alternatively, is it possible that $\mu[S=\text{entire square}]=1$?
Is $S$ a Cantor set $\mu$-a.s.?

\begin{problem}[Gil Kalai, personal communication]
Prove that the limiting measure $\mu$ exists and determine properties of samples
from $\mu$.
\end{problem}

Kalai suspects (personal communication) that
the set $S$ is similar to the set of pivotal sites
(which is a.s.\ a Cantor set in the scaling limit).
This is supported by the easily verified fact that
$\Pb{i\in S}=\Pb{i\text{ pivotal}}$
and $\Pb{i,j\in S}=\Pb{i,j\text{ pivotal}}$
hold for arbitrary boolean functions.
Examples of functions where the scaling limit of $S$
has been determined are provided by Tsirelson~\cite{\TsirelsonFourierWalsh,\TsirelsonScalingLimit}.

\bigskip

One may try to study a scaling limit of dynamical percolation.
Consider dynamical critical site percolation on a triangular lattice of mesh $\eps$, where the
rate at which the sites flip is $\lambda>0$.
We choose $\lambda=\lambda(\eps)$ so that
the correlation between having a left-right crossing of a fixed square
at time $0$ and at time $1$ is $1/2$, say.
Noise sensitivity of percolation~\cite{\BKSnoise} shows that
$\lim_{\eps\searrow 0}\lambda(\eps)=0$ and the results of~\cite{\SSexceptional}
imply that $\lambda(\eps)= \eps^{O(1)}$ and $\eps=\lambda(\eps)^{O(1)}$
for $\eps\in (0,1]$.
It is not hard to invent (several different) notions in which
to take the limit  of dynamical percolation as $\eps\searrow 0$.

\begin{problem}\label{p.dynexists}
Prove that the scaling limit of dynamical critical percolation exists.
Prove that correlations between crossing events at different times $t_1<t_2$
decay to zero as $t_2-t_1\to\infty$ and that a change in a crossing event becomes
unlikely if $t_2-t_1\to 0$.
\end{problem}

Since the correlation between events occuring at different times can be expressed in
terms of the Fourier coefficients~\cite{\BKSnoise,\SSexceptional},
it follows that the second
statement in Problem~\ref{p.dynexists} is very much related to strong concentration
of the measure $\tilde\mu_{f_\ell}$, in the spirit of Problem~\ref{p.muell}.

\medskip

Because of the dependence of $\lambda$ on $\eps$, it is not reasonable to
expect the dynamical percolation scaling limit to be invariant under 
maps of the form $f\times \text{identity}$, where $f:D\to D'$
is conformal and the identity map is applied to the time coordinate.
In particular, in the case where $f(z)=a\,z$, $a>0$,
one should expect dynamical percolation to be invariant under the map
$f\times (t\mapsto a^\beta\,t)$,
where $\beta:= -\lim_{\eps\searrow 0}\log \lambda(\eps)/\log\eps$
(and this limit is expected to exist).
It is not too hard to see that $\beta=3/4$ if the answer to the first question
in Problem~\ref{p.muell} is yes.

This suggests a modified form of conformal invariance for dynamical percolation.
Suppose that $F(z,t)$ has the form $F(z,t)=\bl(f(z),g(z,t)\br)$,
where $f:D\to D'$ is conformal and $g$ satisfies $\pa_t g(z,t)=|f'(z)|^\beta$,
with the above value of $\beta$.
Is dynamical percolation invariant under such maps?
If such invariance is to hold, it would be in a \lq\lq relativistic\rq\rq\ framework,
in which one does not consider crossings occuring at a specific time slice, but rather
inside a space-time set. It is not clear if one can make good sense of that.

\section{LERW and UST}

The loop-erased random walk and the uniform spanning
tree are models where very detailed knowledge
exists. They may be studied using random walks
and electrical network techniques,
and in the two-dimensional setting also by
\SLE/ as well as domino tiling methods.
(See~\cite{\LyonsUSFsurvey} and the references
cited there.)
However, some open problems still remain.

\medskip

One may consider the random walk in a
fine mesh lattice in the unit disk, which
is stopped when it hits the boundary of the disk.
The random walk converges to Brownian motion
while its loop-erasure converges to \SLEk2/.
It is therefore reasonable to expect that the law of the
pair $(\text{random walk},\text{its loop-erasure})$
converges to a coupling of Brownian motion and \SLEk2/.
(If not, a subsequential limit will converge.)

\begin{problem}\label{p.bmsle}
In this coupling, is the \SLEk2/ determined by the
Brownian motion?
\end{problem}

It seems that this question occurred to several researchers
independently, including Wendelin Werner (personal communication).

Of course, one cannot naively loop-erase the Brownian motion
path, because there is no first loop to erase and there are cases
where the erasure of one loop eliminates some of the other
loops.

\bigskip

It is also interesting to try to extend some of the
understanding of probabilistic statistical physics
models beyond the planar setting to higher genus.
The following problem in this direction was
proposed by Russell Lyons (personal communication).

Consider the uniform spanning tree on a fine square grid approximation
of a torus. There is a random graph dual to the tree, which
consists of the dual edges perpendicular to primal edges not in the
tree. It is not hard to see that this random dual of the tree 
contains precisely three edge-simple closed paths
(i.e., no repeating edges), and that these
paths are not null-homotopic.

\begin{problem}\label{p.ustgenus}
Determine the distribution of the triple of homotopy classes
containing these three closed paths.
\end{problem}

The problem would already be interesting for a square torus, but
one could hope to get the answer as a function of the geometry of the torus.

\section{Non-discrete problems}

In this section we mention some problems about the behavior of
\SLE/ itself, which may be stated without relation to
any particular discrete model.

The parametrization of the \SLE/ path by capacity is
very convenient for many calculations. However, in some
situations, for example when you consider the reversal
of the path, this parametrization is not so useful.
It would be great if we had an understanding
of a parametrization by a kind of Hausdorff measure.
Thus we are led to

\begin{problem}\label{p.Hau}
Define a Hausdorff measure on the \SLE/ path
which is $\sigma$-finite.
\end{problem}

We would expect the measure to be a.s.\ finite on
compact subsets of the plane.

That the Hausdorff dimension of the \SLE/ path
is $\min\{2,1+\kappa/8\}$ has been established
 by Vincent Beffara~\cite{\BeffaraSLEdim}.
When $\kappa<8$ (in which case the path has zero area),
we expect the $\sigma$-finite Hausdorff measure to
be the Hausdorff measure with respect to the gauge
function $\phi(r)=r^d\,\log\log(1/r)$, where $d=1+\kappa/8$
is the Hausdorff dimension.
This is based on past experience with similar
random paths~\cite{\TaylorSurvey}.
However, in order to prove that this Hausdorff measure
is $\sigma$-finite, one should probably find alternative
constructions of the measure.
 In the case $\kappa\le 4$,
where the \SLE/ path is a simple path a.s.,
one could try to use conformal maps
from the unit disk to the two components
in the complement of the curve in the upper half plane.
If $f$ is such a map, it might be possible to show
that the limit of the length measure of the
image of the circle $r\,\pa\U$, rescaled appropriately,
has a limit as $r\nearrow 1$. 
Another approach, which was discussed by Tom Kennedy~\cite{\KennedyVariation},
would be to study the $\alpha$-variation of the \SLE/ path,
though this seems hard to handle.

One may also consider other measures of growth for the
\SLE/ path. For example, when $\kappa>4$, we may study the
area of the \SLE/ hull. It would be interesting to study
the various relations between different measures of growth.

\bigskip
It is also natural to ask what kind of sets are visited by
the \SLE/ path. More precisely:

\begin{problem}\label{p.hit}
Fix $\kappa<8$.  Find necessary or sufficient conditions
on a deterministic compact set $K\subset\closure\H$
to satisfy $\Pb{K\cap\gamma\ne\emptyset}>0$,
where $\gamma$ is the \SLEk\kappa/ path.
\end{problem}

The case $K\subset\R$ is of particular interest.
\medskip

When $\kappa=8/3$ and $\H\setminus K$ is
simply connected, there is a simple explicit
formula~\cite{\LSWrestriction} for
$\Pb{\gamma\cap K\ne\emptyset}$.
It is not clear if such formulas are also available
for other values of $\kappa$.
Wendelin Werner~\cite{\WernerLoops} proved the
existence of a random collection of \SLEk{8/3}/-like
loops with some wonderful properties.
In particular, the expected number of loops which
separate two boundary components of an annulus
is conformally invariant, and therefore a function
of the conformal modulus of the annulus.
However, this function is not known explicitly.

\bigskip

Many of the random interfaces which are known or believed
to converge to \SLE/ are reversible, in the sense that the
reversed path has the same law as the original path
(with respect to a slightly modified setup).
This motivates the following problem from~\cite{\RSsle}.

\begin{problem}\label{p.reverse}
Let $\gamma$ be the chordal \SLEk\kappa/ path, where
$\kappa\le 8$.
Prove that up to reparametrization, the image of $\gamma$
under inversion in the unit circle
(that is, the map $z\mapsto 1/\bar z$) has the same law as $\gamma$
itself.
\end{problem}

The reason that we restrict to the case $\kappa\le 8$ is that
this is known to be false as stated when $\kappa>8$~\cite{\RSsle}.
Indeed, there are no known models from physics that are believed
to be related to \SLEk\kappa/ when $\kappa>8$.
Sheffield (personal communication) expects that at least
in the case $\kappa<4$ Problem~\ref{p.reverse} can be 
answered by studying the relationship between the
Gaussian free field and \SLE/.


\bigskip\noindent{\bf Acknowledgments:}
Greg Lawler, Wendelin Werner and Steffen Rohde have collaborated
with me during the early stages of the development of SLE. 
Without them the subject would not be what it is today.
I wish to thank Itai Benjamini, Gil Kalai, Richard Kenyon,
Scott Sheffield, Jeff Steif
and David Wilson for numerous inspiring conversations.
Thanks are also due to Yuval Peres for useful advice, especially concerning
Problem~\ref{p.Hau}.

\bibliography{mr,prep,notmr}
\bibliographystyle{a}

\end{document}